\documentstyle[12pt,amstex]{article}

\makeatletter
\def\rom#1{\leavevmode\skip@\lastskip\unskip\/%
        \ifdim\skip@=\z@\else\hskip\skip@\fi
   {\normalshape#1}}
\makeatother

\newtheorem{Theorem}[equation]{Theorem}

\newcommand{\Hilb}[2]{{#1}^{\lbrack{#2}\rbrack}} 
\newcommand{\Hilbn}[1]{\Hilb{#1}{n}}
\newcommand{\HilbX}[1]{\Hilb{X}{#1}}
\newcommand{\Supp}{\operatorname{Supp}} 

\title{Instantons and affine Lie algebras\footnote{
to appear in Proceedings of Trieste Conference on S-duality and Mirror
Symmetry}}

\author{Hiraku Nakajima\\
Department of Mathematical Sciences,
University of Tokyo, \\
3-8-1 Komaba, Meguro-ku, Tokyo 153, Japan}

\date{September, 1995}

\begin{document}


\maketitle

\begin{abstract}
Various constructions of the affine Lie algebra action on the
moduli space of instantons on $4$-manifolds are discussed.
The analogy between the local-global principle and the role of mass is
also explained.
The detailed proofs are given in separated papers
\cite{Na-algebra,Na-Hilbert}.
\end{abstract}

\section{Introduction}

Vafa and Witten \cite{VW} introduced topological invariants\footnote{
The author does not know how to define their invariants in
a mathematically rigorous way. The difficulty lies in the lack of the
compactness of relevant moduli spaces.}
for $4$-manifolds using $N = 4$ topological supersymmetric Yang-Mills
theory.
Then the $S$-duality conjecture implies that
the generating function of those invariants is a
modular form of certain weight, where the summation runs over
all $\operatorname{SU}(2)$ or $\operatorname{SO}(3)$-principal bundles
of any topological types.
(In general, it has the modular invariance
only for $\Gamma_0(4)$, a subgroup of $\operatorname{SL}(2,\Bbb Z)$.)
For some $4$-manifolds, they identify those invariants with the
Euler numbers of the instanton moduli spaces.
Then they can check the modular invariance
for various $4$-manifolds, using mathematical results,
Yoshioka's formulae \cite{Yo} for ${\Bbb P}^2$ and the blow up,
G\"ottsche and Huybrechts's result \cite{GotHu} for the K3 surface,
and the author's result for ALE spaces.

On the other hand, the author's motivation of the study
\cite{Na-quiver,Na-gauge,Na-algebra} of the homology groups of the
instanton moduli spaces on ALE spaces is totally different.
The author's motivation was trying to understand Ringel \cite{Ri} and
Lusztig's \cite{Lu1,Lu2} constructions of the lower triangular part
${\bold U}_q^-$ of the quantized enveloping algebra.
They used the moduli spaces of representations of quivers, and their
cotangent bundle\footnote{They are not cotangent bundles rigorously.
The situation is very much similar to the relation between the moduli
space of vector bundles over a curve and Hitchin's moduli space of
Higgs bundles.}
can be identified with the instanton moduli spaces on ALE spaces,
via the ADHM description \cite{KN}.
The author showed that the generating function of the Euler numbers of
the instanton moduli spaces on ALE spaces
becomes the character of the affine Lie algebra, which has been
known to have modular transformation property by
Kac-Peterson (see \cite{Kac}).

The definition of the affine Lie algebra representation on the
homology group of the instanton moduli spaces is very geometric, and
seems to be generalized, at least, to projective surfaces.
The results of similar direction are announced recently by
Ginzburg-Kapranov-Vasserot \cite{GKV} and Grojnowski \cite{Gj}.
Unfortunately, our construction depends heavily on the complex
structure of the base manifold.
It is a challenging problem to generalize the construction to more
general $4$-manifolds.
One may need to reformulate the homology group of the moduli spaces, etc.....

\subsection*{Acknowledgements}
The author's understanding of the analogy between the local-global
principle and the theory of mass came from lectures given by the
seminar on Seiberg-Witten theory organized by K.~Ueno.
He would like to speakers, especially S.-K.Yang.
He also thank to G.~Moore and K.~Yoshioka for valuable discussions.

\section{The Hilbert scheme of points and the Heisenberg algebra:
  Twist around points}

In this section, we study the relationship between the Hilbert scheme
of points and the Heisenberg algebra.
The reasons why we study the Hilbert scheme are (a) it is a toy model
for moduli spaces of instantons, (b) it appears in the boundary of
the compactification of the instanton moduli spaces over projective
surfaces, and (c) its homology group is isomorphic to that of a moduli
space for some special cases \cite{GotHu}.

We explain the reason~(b) a little bit more.  Since the instanton
moduli spaces are usually noncompact, one must compactify them to
consider their Euler numbers.  When the base manifold is a projective
surface, the results of Donaldson and Uhlenbeck-Yau enable us to
identify the instanton moduli space with the moduli space of
$\mu$-stable holomorphic vector bundles (Hitchin-Kobayashi
correspondence). Then the one of the most
natural compactifications seems to be Gieseker-Maruyama's
compactifications $\overline{\frak M}$\nobreak
{}~\footnote{The Gieseker-Maruyama's compactifications are not smooth in
  general. In fact, Vafa-Witten's formula for the K3 surface gives the
  fractional Euler number. This may be the contribution of the
  singularities.},
namely moduli spaces of semi-stable torsion free sheaves.

If $\cal E$ is a torsion free sheaf which is not locally free, its
double dual ${\cal E}^{\vee\vee}$ is a locally free sheaf and we have
an exact sequence
\begin{equation*}
  0 @>>> {\cal E} @>>> {\cal E}^{\vee\vee} @>>>
  {\cal E}^{\vee\vee}/{\cal E} @>>> 0.
\end{equation*}
Thus ${\cal E}$ can be determined by (a) ${\cal E}^{\vee\vee}$ and
(b) ${\cal E}^{\vee\vee}\to{\cal E}^{\vee\vee}/{\cal E}$.
The double dual ${\cal E}^{\vee\vee}$ is contained in the interior of
$\overline{\frak M}$, but in the different component with lower second
Chern number.
Thus it is natural to expect that those studies
can be decomposed into two parts, the interior~(a) and
the quotient map~(b).
And the variety of the quotient map~(b),
which depends only on the rank of $\cal E$ and
the length of ${\cal E}^{\vee\vee}/{\cal E}$, looks very much like the
Hilbert scheme of points.
In fact, the Hilbert scheme is the special case
${\cal E}^{\vee\vee} = {\cal O}$.
The Betti numbers of the variety was computed by
Yoshioka~\cite[0.4]{Yo}.

Let $X$ be a projective surface defined over $\Bbb C$.
Let $\Hilbn{X}$ be the component of the
Hilbert scheme of $X$ parameterizing the ideals of ${\cal O}_X$ of
colength $n$.
It is smooth and irreducible.
Let $S^n X$ denotes the $n$-th symmetric product of $X$.
It parameterizes formal linear combinations $\sum n_i [x_i]$ of points
$x_i$ in $X$ with coefficients $n_i\in{\Bbb Z}_{> 0}$ with $\sum n_i = n$.
There is a canonical morphism $\pi\colon \Hilbn{X}\to S^n X$
defined by
\begin{equation*}
  \pi({\cal J}) = \sum_{x\in X}
           \operatorname{length}({\cal O}_X/{\cal J})_x [x].
\end{equation*}
It is known that $\pi$ is a resolution of singularities.

G\"ottsche \cite{Got} computed the generating function of the
Poincar\'e polynomials
\begin{equation*}
  \begin{split}
   & \sum_{n=0}^\infty q^n P_t(\Hilbn{X}) \\
  =  &\prod_{m=1}^\infty
  \prod_{i=0}^4 \, (1 - (-t)^{2m - 2 + i}q^m)^{(-1)^{i+1}b_i(X)}\, ,
  \end{split}
\end{equation*}
where $b_i(X)$ is the Betti number of $X$.
It was shown that the Euler number of $\Hilbn{X}$ is equal to
the orbifold Euler number of $S^n X$ by Hirzebruch-Hofer \cite{HH}.
It was also pointed out by Vafa and Witten that this is equal to the
character of the Fock space.

We shall construct the representation of the Heisenberg and
Clifford algebras in a geometric way. The key point is to introduce
appropriate ``Hecke correspondence'' which give the generators of the
Heisenberg/Clifford algebra.

Take a basis of $H_*(X)$ and
assume that each element is represented by a (real) compact submanifold
$C^a$. ($a$ runs over $1, 2, \dots, \dim H_*(X)$.)
Take a dual basis for $H_*(X)$, and assume
that each element is also represented by a submanifold $D^a$.
(Those assumptions are only for the brevity. The modification to
the case of cycles is clear.)
For each $a = 1,2,\dots,\dim H_*(X)$, $n = 1,2,\dots$ and
$i=1,2,\dots$, we introduce cycles of products of the Hilbert
schemes by
\begin{equation*}
  \begin{split}
    &E_i^a(n) = \{\, ({\cal J}_1,{\cal J}_2)
    \in\HilbX{n-i}\times\HilbX{n} \mid {\cal J}_1\supset {\cal J}_2 \\
    & \quad\text{and
      $\Supp({\cal J}_1/{\cal J}_2) = \{ p\}$ for some $p\in D^a$}\;
      \}, \\
    &F_i^a(n) = \{\, ({\cal J}_1,{\cal J}_2)\in\HilbX{n+i}\times\HilbX{n}
    \mid {\cal J}_1\subset {\cal J}_2 \\
    & \quad\text{and
      $\Supp ({\cal J}_2/{\cal J}_1) = \{ p \}$ for some $p\in C^a$}\;\}.
  \end{split}
\end{equation*}
Then we define an endomorphism $H_*(\HilbX{n}) \to H_*(\HilbX{n-i})$ by
\begin{equation*}
  c\mapsto (p_1)_* ([E_i^a(n)]\cap p_2^* c),
\end{equation*}
where $p_1$, $p_2$ are projections of the first and second factor of
$\HilbX{n-i}\times\HilbX{n}$ and $p_2^* c = [\HilbX{n-i}]\times c$ and
$(p_1)_*$ is a push-forward.
Similarly, we have an endomorphism
$H_*(\HilbX{n}) \to H_*(\HilbX{n+i})$ using $F_i^a(n)$.
Collecting the operators with respect to $n$, we have operators
$[E_i^a]$, $[F_i^a]$ acting on the direct sum
$\bigoplus_n H_*(\HilbX{n})$.
Then
\begin{Theorem}. The following relations hold as operators on
  $\bigoplus_{n} H_*(\HilbX{n})$.
  \begin{gather*}
    [E_i^a] [E_j^b] =
    (-1)^{\dim D^a\dim D^b}[E_j^b] [E_i^a]\\
    [F_i^a] [F_j^b] =
    (-1)^{\dim C^a\dim C^b}[F_j^b] [F_i^a]\\
    [E_i^a] [F_j^b] =
    (-1)^{\dim D^a\dim C^b}[F_j^b] [E_i^a]\\
    \qquad\qquad+ \delta_{ab}\delta_{ij}c_i \operatorname{Id},
  \end{gather*}
  where $c_i$ is a
  nonzero integer depending only on $i$ \rom(independent of $X$\rom).
\label{th:main}\end{Theorem}
In particular, for each fixed $a$, $[E_i^a]$, $[F_i^a]$ ($i =
1,\dots$) define the action of the Heisenberg or Clifford algebra
according to the parity of $\dim C^a$.

Moreover, comparing with G\"ottsche's formula, we can conclude our
representation is irreducible.

The definition of the correspondence $E_i^a$, $F_i^a$ can be naturally
generalized to the case of moduli spaces of torsion free sheaves (see
\cite{Gj}).
However, the author has no idea to generalize to more general
$4$-manifolds.

\section{Elementary transformation: Twist along embedded submanifolds}

The opearator of the previous section twists sheaves around points.
There is another kind of operator which twists sheaves along an
embedded $2$-dimensional submanifold.
This operation is called the elementary transformation in the
literature.

Suppose $C$ is a holomorphic curve embedded in a projective surface
$X$.
Let $i$ denote the inclusion map.
Let $\frak M$ be the moduli space of $\operatorname{U}(r)$-instantons,
namely Einstein-Hermitian connections.
We identify it with the moduli space of holomorhic vector bundles over
$X$ by the Hitchin-Kobayashi correspondence.
It decomposes by the first and second Chern classes.
For each integer $d$, we also consider the following
moduli space ${\frak P}$ of parabolic bundles
$({\cal E_1},{\cal E_2},\varphi)$ where
$\cal E_i$ is a holomorphic vector bundle over $X$, and
$\varphi\colon{\cal E_1} \to {\cal E_2}$ is an injection which is an
isomorphism outside $C$.
In order to define the moduli space, we need to introduce the notion
of the stability to parabolic bundles (see \cite{MY}),
for this we need to choose an ample line bundle $L$ or the K\"ahler
metric which is a curvature of $L$.
Moreover, it is necessary to consider the Gieseker-Maruyama
compactification of moduli spaces, as explained in the previous
section.
But we do not go in detail.

There is a morphism $f\colon {\frak P}\to {\frak M}\times {\frak M}$.
Then we can define two operators on the homology group exactly as in the
previous section:
\begin{equation*}
\begin{split}
  c &\mapsto (p_1)_* (f_*[{\frak P}]\cap p_2^* c),\\
  c &\mapsto (p_2)_* (p_1^* c\cap f_*[{\frak P}])
\end{split}
\end{equation*}
Note that the first and second Chern classes are not preserved this
operator.
Strictly speaking, we do not have a globally defined morphism
since the stability conditions for the parabolic bundles
and their underlying vector bundles are not equivalent in general.
But it is enough for our purpose to have $f_*[{\frak P}]$ as an
element of homology group of ${\frak M}\times {\frak M}$.
For example, $f$ could be a meromorphic map.

Since we do not know what is the right setting for general
projective surfaces\footnote{If one could define the Hecke operators using
  Kronheimer-Mrowka's singular anti-self-dual connections, they might be
  the right setting.}, we focus on particular examples, namely ALE
spaces.
The ALE spaces are the minimal resolution of simple singularities
${\Bbb C}^2/\Gamma$, where $\Gamma$ is a finite subgroup of
$\operatorname{SU}(2)$.
The second homology group $H_2$ of the ALE space is spanned by
the irreducible components $\Sigma_1,\dots, \Sigma_n$ of the exceptional set,
which are the projective line.
The intersection matrix is the negative of the Cartan matrix of type
ADE.
The classification of simple singularities are given by the Dynkin
graphs in this way.
In particular, there is a bijection between simple singularities and
simple Lie algebra of type ADE.
The rank is equal to the number of the number of the irreducible
components, namely $n$.
By the work of Kronheimer, it is known that
they have hyper-K\"ahler metrics.
There are a variant of the ADHM description, which identifies the
framed moduli spaces of instantons, or more precisely, torsion-free
sheaves with the cotangent bundles of the moduli space
of representations of quivers of affine Dynkin graphs.

Since we shall work on non-compact spaces, we have extra discrete
parameters which parameterizes the boundary condition.
We consider instantons which converge to a flat connection at the end
of the ALE space $X$.
The flat connection on the end can be classified by its monodromy,
namely a representation $\rho$ of the finite group $\Gamma$.
Let $\rho_0$, $\rho_1$, \dots, $\rho_n$ be the irreducible
representations of $\Gamma$ with $\rho_0$ the trivial representation.
By the McKay correspondence, there is a bijection between the vertices
of the affine Dynkin graph and the irreducible representations
(see \cite{Na-gauge} for more detail).
The monodromy representation $\rho$ is decomposed as
$\rho = \bigoplus_{k=0}^n \rho_k^{\oplus w_k}$, where $w_k$ is the
multiplicity. This datum will be preserved under the Hecke operator.

Corresponding to each irreducible component $\Sigma_k$,
we take a component of the moduli space of parabolic bundles where
${\cal E_2}/\varphi({\cal E_1})$ is rank $1$ and degree $-1$.
We then define operators
$e_k$ and $f_k$ on the homology group of the moduli space as above.
We also have an operator $\alpha_k^\vee$ which is the multiplication
by $-\langle c_1, [\Sigma_k]\rangle$ on the homology class belonging
to the component with the first Chern class $c_1$.

For $k = 0$, we can define similar operators $e_0$, $f_0$ by replacing
${\cal O}_{\Sigma_k}(-1)$ by a sheaf ${\cal O}_{\bigcup \Sigma_k}$.
The operator $\alpha_0^\vee$ is defined so that
\begin{equation*}
   \sum_{k=0}^n \dim\rho_k \alpha_k^\vee
   = \operatorname{rank} E\, \operatorname{Id}.
\end{equation*}

Finally define the operator $d$ to detect the instanton number.
Namely the mulitiplication by
\begin{equation*}
  -\int_X \operatorname{ch}({\cal E}).
\end{equation*}
on the homology group of the each component of the moduli spaces.

\begin{Theorem}.
Operators $\alpha_k^\vee$, $e_k$, $f_k$ \rom($k = 0, \dots, n$\rom), $d$
satisfy the relation of the affine Lie algebra corresponding to the
extended Dynkin graph.
Moreover, the representation on $H_*({\frak M})$ is integrable.
\end{Theorem}

The irreducible decomposition of the representation is complicated,
but we have one irreducible factor whose geometric meaning is clear.

\begin{Theorem}.
If we take the middle degree part of the homology group
$H_*({\frak M})$ \rom(since the dimension of $\frak M$ are changing on
components, the middle degree also changes\rom), it is preserved by
the affine Lie algebra action.
Moreover, it is the integrable highest weight representation with the
highest weight vector $\,{}^{t}(w_0,\dots,w_n)$.
The level is equal to the rank of the vector bundle.
\end{Theorem}

The highest weight vector lies in the paticularly chosen moduli space
which consists of a single point.

\section{Local-Global Principle and the Mass}

Historically the Hecke operators were originally introduced in the
theory of modular forms. There are also analogoues operators in the
theory of the moduli spaces vector bundles over curves \cite{NR},
which are used in the geometric Langlands program.
Our operators can be considered as natural complex $2$-dimensional
analogue of these operators.

The importance of the Hecke operators comes from the fact
that they lie in the heart of the ``local-global principle''.
We shall explain it only very briefly.
The interested reader should consult to good literatures about the
Hecke operators and the modular forms (see e.g., \cite{Langlands}).

The local-global principle roughly says that
a global problem could be studied as a collection of local
problems.
The basic example is the Hasse principle: A quadratic from with
integer coefficients has a nontrivial integer solution if and only if it
has real solution and a $p$-adic solution.
In this case the global problem is to find an integer solution and the
local problems are find solutions in $\Bbb R$ and $\Bbb Q_p$.
Thus the base manifold, which is the parameter space of the local
places, is the set of prime numbers plus infinity $\Bbb R$.

The theory of modular forms are also examples of the global-local
principle.
Consider the space of modular forms of weight $k$, which can be
considered as functions of lattices in
$\Bbb C$ with homogeneous degree $-k$, i.e., $F(\lambda L) =
\lambda^{-k}F(L)$.
Then for each prime number $p$, we define the Hecke operator $T(p)$ by
\begin{equation*}
  (T(p) F)(L) = p^{k-1} \sum_{[ L' : L] = p} F(L'),
\end{equation*}
where the summation runs over the set of sublattices of $L$ with index
$p$.
These operators commute each other.
If a modular form is a simultaneous eigenfunction, its $L$-function
has an Euler product expansion. The analogy between the modular forms
and $4$-dimensional gauge theory our theory are given in the table.

\begin{table}[hbt]
\begin{tabular}{l|l}
\hline
prime numbers        & points in a $4$-manifold \\
                     & \ and submanifolds $C$\\
\hline
the space of         & the homology gruop \\
\ modular forms      & \ of moduli spaces \\
\hline
Hecke operators      & our Hecke operators\\
\hline
\end{tabular}
\end{table}

In physics, there is a good explanation why the local-global principle
holds in some topological field theories\footnote{Topological field
  theories of cohomological type according to the terminology in
  \cite{CMR}}. In these theories, topological invariants, like
Donaldson's invariants, are expressed as correlation functions.  It is
not by no means obvious that the results are topological invariants,
since one needs to introduce a Riemannian metric for the definition of
the Lagrangian.  However, by a clever choice of the Lagrangian, the
resulted correlation functions are independent of the choice of the
metric. The mechanism is just like the fact that the euler class,
which is defined as the pfaffian of the curvature, is a topological
invariant.  Thus one can take a family of Riemannian metrics $g_t =
tg$ with $t > 0$, and consider the limiting behaviour for $t\to
\infty$. In the limit, the distance of two different points goes to
infinity. Hence if the ``mass'' of all particles is not zero, there
are no interaction between two points. Then one can compute the
correlation functions by integrating local contributions over the base
manifold. In this sense, the local-global principle holds in this
theory.

In $N = 1$ topological supersymmetric Yang-Mills theory, it is
believed that all particles have non-trivial mass. However in $N = 2$
topological supersymmetric Yang-Mills theory, which is relevant to
Donaldson's invariants, it is no longer true. Hence there may be
massless particles which make interaction even when the distance
between two points are very large.  Anyway, if all particles would
have non-trivial mass,
Donaldson's invariants would depend only on homology classes of the
underlying manifolds.  On a K\"ahler manifold $X$
with a non-trivial holomorphic $2$-form $\omega$, Witten
\cite{Wi-SUSY} used a
perturabation from the $N = 2$ theory to the $N = 1$ theory adding a
term depending on $\omega$. The remarkable observation was that there
remain particles which have zero mass where $\omega$ vanishes.
Unless the manifold $X$ is a K3 surface or a
torus, $\omega$ vanishes along a divisor $C$.
Thus the local-global principle holds in the $N = 2$ theory with only
one modification; there are non-trivial contribution from $C$.
In other words, the $2$-dimensional submanifold $C$ cannot be divided
any more, and should be considered as a point.
Similarly, in \cite{VW}, again for the same class of K\"ahler
manifolds, the $N = 4$ theory was perturbed to the $N = 1$ theory, and
the correlation function was calculated in a similar way.
And finally, Witten conjectured \cite{Wi-monopole}
that the Kronheimer-Mrowka's basic classes \cite{KrMr}
coincide with homology classes whose Seiberg-Witten invariants are
nonzero.
It means that the local-global principle fails exactly along basic
classes.

Now it becomes clear why we must introduce two kinds of Hecke
correspondences, twist along a point and twist along a $2$-dimensional
submanifold.
For the Hilbert scheme on a projective surface, it is apparent that
the local-global principle holds without the introduction of $C$.
This should be the basic reason why the homology group of the Hilbert
scheme is generated only by the first kind of the Hecke correspondence.
For higher rank case, we mighty need the second type of the
Hecke correspondence in order to get all homology classes in the
moduli spaces.
However, the relation between two kinds of Hecke operators is not clear, at
this moment.
Moreover, there might be other types of correspondences which is
useful to describe the homology group of the moduli space.
For example, G\"ottsche and Huybrechts
\cite{GotHu} used an interesting correspondence in order to relate moduli
spaces of rank $2$ bundles and Hilbert schemes.

Finally, we would like to point out the difference between our
situation and the classical one (i.e., the Hecke operators on
modular forms).
The first kind of the Hecke
operator is independent of the choice of the representative $C^a$ of
the homology class.
This is because we are studing the homology group of the moduli space.
It is not clear that the second kind of the Hecke operator depends
only on the homology class of $C$. But it depends only on the rational
equivalence class.

\end{document}